\documentclass{article}
\usepackage{kluwerstatandcomp}
\usepackage{amsmath, amssymb}
\usepackage[hypertex]{hyperref}
\usepackage{natbib}
\usepackage{graphicx,   epic, eepicemu, epsfig, float, enumerate}
\newtheorem{theorem}{Theorem}[section]
\newtheorem{example}[theorem]{Example}
\newtheorem{remark}[theorem]{Remark}
\def\V{\mathcal{V}}

\bibpunct{(}{)}{,}{a}{}{,}

\begin{document}\allowdisplaybreaks
\title{{\huge \bf Adjusted Viterbi Training}}
\author{{\Large J\"uri Lember}\thanks{Estonian Science Foundation Grant 5694}\\
Tartu University, Liivi 2-507, Tartu 50409, Estonia; jyril@ut.ee\\
{\Large Alexey Koloydenko}\thanks{Corresponding author}\\
Eurandom, P.O. Box 513 - 5600 MB Eindhoven, The Netherlands.\\
Tel: +31(0)40.247.8129, koloidenko@eurandom.tue.nl\\
\href{http://euridice.tue.nl/~akoloide/VA/}{http://euridice.tue.nl/$\sim$akoloide/VA/}}
\date{January 20, 2004}
\maketitle
\begin{abstract}
We study modifications of the Viterbi Training (VT) algorithm to estimate emission 
parameters in Hidden Markov Models (HMM) in general, and in mixure models in particular. 
Motivated by applications of VT to HMM that are used in speech recognition, natural 
language modeling,  
image analysis, and bioinformatics, we investigate a possibility of alleviating the 
inconsistency of VT while controlling the amount of extra computations.  Specifically, 
we propose to enable VT to 
asymptotically fix the true values of the parameters as does the EM algorithm.  
This relies on infinite Viterbi alignment and an associated with it limiting probability 
distribution. This paper, however, focuses on mixture models, an important 
case of HMM, wherein the limiting distribution can always be computed exactly;
finding such limiting distribution for general HMM presents a more challenging
task under our ongoing investigation.

A simulation of a univariate Gaussian mixture shows that our central algorithm (VA1) 
can dramatically improve accuracy without much cost in computation time.

We also present VA2, a more mathematically advanced correction to VT, verify by simulation
its fast convergence and high accuracy; its computational feasibility remains to
be investigated in future work.
\end{abstract}
\section{Introduction}
\label{sec:intro}
Motivated by applications of the Viterbi Training (VT) algorithm to 
estimate parameters of Hidden Markov Models in speech  recognition 
\citep{vanaraamat, philips,  raamat, Rabiner86, philips2, strom},
natural language modeling \citep{ochney},  image analysis \citep{gray2}, bioinformatics 
\citep{dna2, dna1}, and  by connections with constrained vector quantization \citep{gray, gray4},
we are interested in improving accuracy of the VT estimators while preserving
its essential computational advantages. 
\\ \indent
Let $\theta_l$ be the emission parameters of HMM with states $l\in S=\{1,...,K\}$.
The standard method to compute a maximum likelihood estimator of  
$\theta$ is the  EM algorithm that in the HMM context is also known as the {\em Baum-Welch} or
{\em forward-backward algorithm} \citep{baumhmm, emtutorial,  vanaraamat,  jelinek, tutorial, raamat,
young_HMM}.
Since EM  can in practice be  computationally expensive, it is commonly replaced  by
{\em Viterbi Training}. VT effectively replaces the computationally costly 
expectation (E) step of  EM  by an appropriate maximization step that is 
computationally less intensive. An important example of successful and 
elaborate application of VT in industry is Philips speech recognition systems 
\citep{philips}. 
\\ \indent There are also variations of VT that use more than one 
best alignment, or several perturbations of the best alignment \citep{ochney}.  
The improvements that we explore are, however, of a different nature: 
roughly, we increase estimation accuracy purely by means of analytic calculations and do 
not require computing more than one optimal alignment. 
\\ \indent
VT is inferior to  EM in terms of accuracy because
the  VT estimators need not be (local) maximum likelihood estimators 
(VT does not necessarily increase the likelihood); 
which leads to bias and inconsistency (\S\ref{sec:viterbi_alignment_training}).  
\\ \indent
Given current parameter values, VT first finds a Viterbi {\em alignment} that is 
a sequence of hidden states maximizing 
the likelihood of the  observed data.  Observations  assumed to have been emitted
from  state $l$, are
regarded as an {\em i.i.d.} sample  from $P_l$, the corresponding emission distribution.
These observations produce $\hat P_l^n$, the empirical version of $P_l$, and ultimately 
$\hat{\mu}_l$, a maximum likelihood estimate of $\theta_l$.
$\hat{\mu}$ is then used to find an alignment in the next step, 
and so forth.  It can be shown that in general this procedure converges in finitely 
many steps; also, it is usually much faster than EM. 
\\ \indent
In speech recognition,  the same training  procedure was already  
described by L. Rabiner {\em et al.} in \citep{Rabiner90, Rabiner86} 
(see also \citep{tutorial, raamat})
who considered his procedure as a variation of the {\em Lloyd algorithm} from vector
quantization, and referred to it as {\em segmential K-means  training}. 
The  analogy  with  vector quantization is  especially pronounced  when the
underlying  chain is a sequence of {\em i.i.d.} variables
in which case the observations  are simply an {\em i.i.d.} sample from a mixture
distribution (\S\ref{sec:mixture}). For such mixture models, Viterbi training  was also 
described  by R. Gray {\it et al.} in \citep{gray}, where the training  algorithm was 
considered in the vector quantization context under the name of {\em entropy constrained
vector  quantization (ECVQ)}.  (See also \citep{gray4} for more recent developments in this theory.)
\\ \indent
Our goal is to alleviate the inconsistency of the VT estimators while
preserving the fast convergence and 
computational feasibility of the baseline VT algorithm.  To this end, 
we note that $\theta^*$, the true parameters, are asymptotically a fixed point 
of  EM  but not of VT \S\ref{sec:viterbi_alignment_training}, \S\ref{sec:mixture}. 
We thus attempt to adjust VT in order to restore this property by studying the asymptotics 
of $\hat P^n_l$.  Thus,  
we discuss the existence of $Q_l$  $l\in S$
\begin{equation}\label{ujuI}
\hat P^n_l\Rightarrow Q_l,\quad l\in S \quad{\rm a.s.},
\end{equation}
first in the general HMM context -- \S\ref{sec:viterbi_alignment_training}, and then in
the special case of mixture models -- \S\ref{sec:mixture}.
If such limiting measures exist, then under certain continuity assumptions, the estimators
$\hat{\mu}_l$ will converge to $\mu_l$, where
\begin{equation*}\label{aziz}
\mu_l=\arg\max_{\theta_l}\int \ln f_l(\theta_l,x)Q_l(dx).
\end{equation*}
Taking  into account the difference between $\mu_l$ and the true parameter, the appropriate
adjustment of the Viterbi training can now be defined (\S\ref{sec:viterbi_alignment_training}).
\\ \indent
However, the  asymptotic behavior of $\hat P^n_l$ is not in general straightforward
and its analysis requires an extension of the definition of Viterbi alignment 
{\em at infinitum} \citep{lembereur}. With the
infinite alignment one can prove 
the existence of the limiting measures $Q_l$ \eqref{ujuI}, which is essential for the general 
definition of the adjusted Viterbi training.
\\ \indent
To implement these ideas in practice, a closed form of $Q_l$  
(or ${\hat \mu}_l$) as a function of the true parameters is necessary.  The measures $Q_l$ 
also depend on the transition as well as on the emission models, and  computing $Q_l$ 
can be very difficult.  
However, in the special case of mixture models \S\ref{sec:mixture}, the measures $Q_l$ 
are easier to find.  We are also motivated by the continuing interest of others 
in computational efficiency and accuracy of parameter estimation in mixture 
models \citep{dias, LinChenWu}. In \S\ref{sec:mixture}, 
we describe the adjusted Viterbi training (VA1) for the mixture case, which we 
view as the  main contribution of this paper: VA1 {\em recovers the asymptotic fixed point 
property} and, since its adjustment function {\em does not depend on data}, each iteration of 
VA1 enjoys the {\em same order of computational complexity (in terms of the sample size) as VT}. 
Moreover, for commonly used mixtures, such as, for example (Example \ref{normal}), mixtures of 
multivariate normal distributions with unknown means and known covariances, the adjustment function 
is available in a {\em closed form} requiring integration with the mixture densities.  Depending 
on the dimension of the emission variates and on the number of components, and on the available 
computational resources, one can vary the accuracy of  the adjustment.  We reiterate that,
unlike  the computations of the EM algorithm, computations of the adjustment 
{\em do not involve evaluation and subsequent summation of the mixture density at every data point}.
\\ \indent
We first introduce these ideas for the case of known mixture weights 
and then extend them in \S\ref{sec:kaalud} to the case of unknown weights.  
\\ \indent
To test our theory, in \S\ref{sec:simu} we simulate a mixture of two univariate normal 
distributions with unit variance, unknown means, and unequal but comparable weights.  
The main goal of our simulations is  to compare the
performances of VT, VA1, and EM in terms of the  accuracy,  convergence, 
amount of computations per iteration, and the  total amount of computations.
The simulations are performed under different conditions on the initial guess 
when the weights are assumed to be known, \S\ref{sec:known_weights}, and unknown, 
\S\ref{sec:unknown_weights}, and the results (\S\ref{sec:results}) 
are consistently in favor of VA1. 
\\ \indent
In \S\ref{sec:second}, we  propose VA2, a more mathematically advanced correction to VT; 
we verify its fast convergence and high accuracy on the simulated data in \S\ref{sec:simu}, 
and intend to elaborate on its computationally feasible implementations in future work.  
A concluding summary is given in \S\ref{sec:conclusion}.
\section{General HMM}\label{sec:viterbi_alignment_training}
Let $Y$ be a Markov chain with finite state space $S$. We assume  $Y$ to be
 irreducible and aperiodic with transition matrix
$P=(p_{ij})$ and  initial distribution $\pi$ that is also the stationary distribution of $Y$.
 To every state $l\in S$ there corresponds an {\it emission distribution} $P_l$ on
(${\mathcal X}, {\mathcal B}$),  a separable metric
space and the corresponding Borel $\sigma$-algebra. Let $f_l$, the density of $P_l$ 
with respect to some reference measure $\lambda$ (for instance, the Lebesgue measure), be 
known up to the parametrization $f_l(x; \theta_l)$. When $Y$ is in state
$l$, an observation according to $P_l(\theta^*)$ and independent of everything else is emitted,
with $\theta^*=(\theta^*_1,\ldots,\theta^*_K)$ being the unknown true parameters. 
\\ \indent
Thus, for any $y=y_1,y_2,\ldots $, a realization of $Y$,
there corresponds a sequence of independent random variables, $X_1,X_2,\ldots $,
where $X_n$ has  distribution $P_{y_n}$.
Note that we only observe  $X=X_1,X_2,\ldots$ and the realization $y$ is unknown 
($Y$ is hidden).
\\ \indent
The distribution of $X$ is completely determined by
the chain parameters $(P,\pi)$ and the emission distributions
$P_l,$ $l\in S$. The process $X$ is also {\em mixing} and, therefore, ergodic.
We now recall Viterbi Alignment and Training.
\\ \indent
Let $x_1,\ldots,x_n$ be first $n$ observations on $X$. 
Let $\Lambda(q_1,\ldots,q_n;x_1,\ldots,x_n; \theta)$  be the likelihood function
${\bf P}(Y_i=q_i,~i=1,\ldots,n)\prod_{i=1}^nf_{q_i}(x_i;\theta_{q_i})$, $q_i\in S$.
\\ \indent
The {\em Viterbi alignment} is any sequence of states
$q_1,\ldots,q_n\in S$ that maximizes the likelihood of observing $x_1,\ldots,x_n$, 
$\theta$ being fixed.
Thus, for a fixed $\theta$, the Viterbi
alignment is a maximum-likelihood estimator  of {\em the realization of}
$Y_1,\ldots, Y_n$ given $x_1,\ldots,x_n$. In the following, the
Viterbi alignment will be referred to as the alignment. For each  $n\geq 1$, 
let $\V$ be the set of  alignments:
\begin{equation}\label{alignment}
  \V(x_1,\ldots,x_n;\theta)= \{v\in S^n:~
    \forall w\in S^n~\Lambda(v;x_1,\ldots,x_n;\theta)\ge  \Lambda(w;x_1,\ldots,x_n;\theta)\}.
\end{equation}
Any map $v:{\mathcal X}^n\mapsto \V(x_1,\ldots,x_n;\theta)$ will also be called  an alignment. 
Further, unless explicitly specified, $v_{\theta}$ will denote an arbitrary element of 
$\V(x_1,\ldots,x_n;\theta)$.
\begin{center}\underline{Viterbi Training}\end{center}
\begin{enumerate}[1.)]
\item Choose an initial value $\theta^0=(\theta^0_1,\ldots,\theta^0_K)$.
\item Given $\theta^j$ $j\ge0$, find the alignment $$v_{\theta^j}(x_1,\ldots,x_n) =(v_1,\ldots,v_n)$$
and partition  $x_1,\ldots, x_n$ into $K$ subsamples, with $x_k$ going 
to the $l^{th}$ subsample if and only if $v_k=l$. Equivalently,  define  $K$  empirical measures
\begin{equation}\label{emp}
\hat P_l^n(A;\theta^j):={\sum_{i=1}^n I_{A\times l}(x_i,v_i)\over \sum_{i=1}^n I_l(v_i)},
\quad A\in {\mathcal B},\quad l\in S,
\end{equation} where $I_A$ stands for the indicator function of set $A$.
\item For every subsample find the MLE given by:
\begin{equation}\label{mle}
\hat{\mu_l}(\theta^j)=\arg\max_{\theta_l\in \Theta_l} \int \ln f_l(x;\theta_l)\hat P^n_l(dx;\theta^j),
\end{equation}
and take $\theta^{j+1}_l=\hat{\mu_l}(\theta^j),\quad l\in S$.
If for some $l\in S$,  $v_i\ne l$ for any $i=1,\ldots, n$ ($l^{th}$ subsample is empty), then 
the empirical measure $\hat P_l^n$ is formally undefined, in which  case we take $\theta_l^{j+1}=\theta^j_l$.
We omit this exceptional case in the following discussion.
\end{enumerate} 
VT can be interpreted as follows. Suppose that
at step $j$, $\theta^j=\theta^*$ and hence $v_{\theta^j}$ is obtained using the
true parameters. The training is then  based on the assumption that
the alignment  $v(x_1,\ldots,x_n)=(v_1,\ldots,v_n)$ is correct, 
i.e., $v_i=Y_i$, $i=1,\ldots,n$. In this case, the empirical measures $\hat P_l^n(\theta^j)$,
$l\in S$ would be obtained from the i.i.d. sample generated from $P_l(\theta^*)$,  and
the MLE  $\hat{\mu_l}(\theta^*)$ would be the natural estimator to use.
Clearly, under these assumptions $\hat P^n_l(\theta^*)\Rightarrow P_l(\theta^*)$ a.s.
and, provided that $\{f_l(\cdot; \theta):\theta\in
\Theta_l\}$ is a $P_l$-Glivenko-Cantelli class and $\Theta_l$ is
equipped with some suitable metric, $\lim_{n\to\infty}\hat \mu_l(\theta^*)= \theta^*_l$
a.s.  Hence, if $n$ is sufficiently large, then  $\hat P_l^n\approx P_l$ and
$\theta_l^{j+1}=\hat \mu_l(\theta^*) \approx \theta^*_l=\theta^j_l$,
$\forall l$ i.e. $\theta^j=\theta^*$ would be (approximately) a fixed point of the training algorithm.
\\ \indent
A weak point of the previous argument is that the alignment in general is
not correct even when the parameters used to find it are, 
i.e. generally $v_i\ne Y_i$. In particular, this implies
that the empirical measures $\hat P_l^n(\theta^*)$ are not obtained from an
i.i.d. sample taken from $P_l(\theta^*)$.   Hence, we have no reason to
believe that $\hat P_l^n(\theta^*)\Rightarrow P_l(\theta^*)$ a.s. and $\lim_{n\to\infty}\hat{\mu_l}(\theta^*)=
\theta^*_l$ a.s. Moreover, we do not even know whether
the sequences of empirical measures $\{\hat P_l^n(\theta^*\}$ and MLE estimators
$\{\hat{\mu_l}(\theta^*)\}$ converge (a.s.) at all.
\\ \indent
In \citep{lembereur}, we prove the existence of limiting probability measures
$Q_l(\theta,\theta^*)$, $l\in S$, that depend on $\theta$, the parameters
used to find the alignment $v_\theta(x_1,\ldots,x_n)$, and on $\theta^*$,
the true parameters with which the random samples are generated.  Namely,
these $Q_l$, $l\in S$ will be such that for every $l$
\begin{equation}\label{koondumineI}
\hat P_l^n(\theta^*)\Rightarrow Q_l(\theta^*,\theta^*),\quad \text{a.s.}.
\end{equation}
Suppose also that the parameter space $\Theta_l$ is equipped with some metric. Then, under certain
consistency assumptions on classes  ${\mathcal F}_l=\{f_l(\theta_l):\theta_l\in \Theta_l\}$,  
the convergence
\begin{equation}\label{koondumineII}
\lim_{n\to\infty}{\hat \mu}_l(\theta^*) =\mu_l(\theta^*,\theta^*)\quad \text{a.s.}
\end{equation}
can be deduced from \eqref{koondumineI}, where
\begin{equation}\label{koondumineII2}
\mu_l(\theta,\theta^*)\stackrel{\rm def}{=}\arg\max_{\theta'_l\in \Theta_l}\int 
  \ln f_l(x; \theta'_l)Q_l(dx;\theta,\theta^*).
\end{equation}
We also show that  in general, for the baseline Viterbi training 
$Q_l(\theta^*,\theta^*)\ne P_l(\theta^*)$, 
implying  $\mu_l(\theta^*,\theta^*)\ne \theta^*_l$. In an attempt to 
reduce the bias $\theta_l^*-\mu_l(\theta^*,\theta^*)$,  we next propose the 
{\it adjusted Viterbi training}.
Suppose \eqref{koondumineI} and \eqref{koondumineII} hold.
Based on \eqref{koondumineII2}, we now consider the mapping
\begin{equation}\label{mapping}
\mu_l(\theta)=\mu_l(\theta,\theta),\quad l=1,\ldots, K.
\end{equation}
Since this function is independent of
the sample,  we  can define the following correction for the bias:
\begin{equation}\label{mapII}
\Delta_l(\theta)=\theta_l-\mu_l(\theta),\quad l=1,\ldots, K.
\end{equation}
\begin{center}\underline{VA1 -- Adjusted Viterbi Training}\end{center}
\begin{enumerate}[1.)]
\item Choose an initial value $\theta^0=(\theta^0_1,\ldots,\theta^0_K)$.
\item Given $\theta^j$, perform the alignment and define $K$ empirical measures $\hat P_l^n(\theta^j)$ 
as in \eqref{emp}.
\item For every $\hat P_l^n$, find $\hat{\mu_l}(\theta^j)$ as in \eqref{mle} and for each $l$, define
$\theta^{j+1}_l=\hat{\mu_l}(\theta^j)+\Delta_l(\theta^j)$, where $\Delta_l$ is defined $\forall l\in S$
in \eqref{mapII}.
\end{enumerate}
Note that, as desired, for $n$ sufficiently large, the adjusted training algorithm  has $\theta^*$ 
as its (approximately) fixed point: Indeed, suppose $\theta^j=\theta^*$. From \eqref{koondumineII}, 
$\hat{\mu}_l(\theta^j)=\hat{\mu}_l(\theta^*)\approx \mu_l(\theta^*)=\mu_l(\theta^j)$, for all $l\in S$.
Hence, 
\begin{equation}
  \label{eq:approx}
\theta_l^{j+1}=\hat{\mu}_l(\theta^*)+\Delta_l(\theta^*) \approx 
\mu_l(\theta^*,\theta^*)+\Delta_l(\theta^*)=
\theta_l^*=\theta^j,\quad l\in S.  
\end{equation}
\section{Mixture model}
\label{sec:mixture}
In general, no closed form for the distribution
$Q_l(\theta^*,\theta^*)$  in \eqref{koondumineI} is available. Therefore, the
mapping  \eqref{mapping}  may be impossible to determine exactly and 
approximations of $Q_l$ should  be used for the adjustments of Viterbi training 
(\S\ref{sec:viterbi_alignment_training}).
However, in the case of the mixture model, the distributions $Q_l$ are straightforward 
to find and the adjusted Viterbi training can therefore be fully specified.
In this model, $Y$, the underlying Markov chain, is a
sequence of i.i.d. discrete random variables with the state space 
$S=\{1,\ldots,K\}$ of {\em mixture components}.  Thus, the transition probabilities are
$p_{ij}=p_j$, $i,j\in S$, where $p_j$ are mixture weights. 
To each component $l\in S$ there corresponds a probability distribution $P_l$ 
with density $f_l=f_l(\cdot;\theta^*_l)$, where $\theta^*_i$ are the true parameters.
Unless explicitly stated otherwise, 
the mixture weights $p_l$ will be assumed to be known. Such a model produces
observations $x_1,\ldots,x_n$, that are regarded as
an i.i.d. sample from the mixture distribution $P$ with density
\begin{equation}\label{segutih}
\sum_{i=1}^Kp_if_i=\sum_{i=1}^Kp_if_i(\cdot;\theta^*_i)=f(\cdot;\theta^*)=f.
\end{equation}
For any set of parameters $\theta=(\theta_1,\ldots \theta_K)$, the alignment $v_{\theta}$
can be obtained via a {\it Voronoi partition} 
${\mathcal S}(\theta)=\{S_1(\theta),\ldots,S_K(\theta)\}$, where
\begin{align}
S_1(\theta)&=\{x: p_1f_1(x;\theta_1)\geq p_jf_j(x;\theta_j),\quad \forall j\in S\} \label{eq:S1}\\
S_l(\theta)&=\{x: p_lf_l(x;\theta_l)\geq p_jf_j(x;\theta_j),\quad \forall j\in S\}\backslash 
(S_1\cup\ldots\cup S_{l-1}),\quad l=2,\ldots, K.\label{eq:S}
\end{align}
Now, the alignment can be defined as follows: $v_{\theta}(x)=l$ if
and only if $x\in S_l(\theta)$. 
In particular, given  the Voronoi partition ${\mathcal S}(\theta)=\{S_1,\ldots,S_l\}$, 
the empirical measures $\hat P_l^n$ \eqref{emp}  are
\begin{equation}\label{empmix}
\hat P_l^n(A;\theta)={\sum_{i=1}^n I_{S_l(\theta)\cap A}(x_i)\over \sum_{i=1}^n I_{S_l(\theta)}(x_i)},
\quad A\in {\mathcal B},\quad l\in S.
\end{equation}
Thus, given the same partition,
$\hat{\mu}_l(\theta)$ \eqref{mle}, the sub-sample MLE for component $l$, becomes
\begin{equation}\label{hyljes}
\hat{\mu}_l(\theta)=\arg\max_{\theta'_l\in \Theta_l}\int_{S_l(\theta)}\ln f_l(x;\theta'_l)\hat P_n(dx),
\end{equation}
where $\hat P_n$ is the ordinary empirical measure associated with the given random sample.
The convergence \eqref{koondumineI} then  follows immediately from \eqref{empmix}.
Indeed, for any $\theta$, 
by virtue of the {\em Strong Law of Large Numbers}  we have 
\begin{align*}
\lim\limits_{n\to\infty}\hat P_l^n(A;\theta)&\stackrel{\rm a.s.}{=} 
{P(A\cap S_l(\theta);\theta^*)\over P(S_l(\theta);\theta^*)}=
{\int_{S_l(\theta)\cap A}f(x;\theta^*)d\lambda(x)\over \int_{S_l(\theta)}f(x;\theta^*)d\lambda(x)}
=&{\sum_ip_i\int_{S_l(\theta)\cap A}f_i(x;\theta^*_i)d\lambda(x)
\over\sum_ip_i\int_{S_l(\theta)}f_i(x;\theta^*_i)d\lambda(x)}.  
\end{align*}
Since ${\mathcal X}$ is separable, 
it follows that $\hat P_l^n\Rightarrow Q_l$ a.s., where 
$$q_l(x;\theta,\theta^*)\propto f(x;\theta^*)I_{S_l(\theta)}=
(\sum_ip_if_i(x;\theta^*))I_{S_l(\theta)},\quad l=1,\ldots, K$$
are the densities of respective $Q_l(\theta,\theta^*)$'s.  
\\ \indent
Now it is clear that even when the partition ${\mathcal S(\theta^*)}$ is obtained using 
the true parameters $\theta^*$, 
$Q_l(\theta^*,\theta^*)$, the limiting distribution (density $q_l(x;\theta^*,\theta^*)$), 
can be different from $P_l(\theta^*)$, the desired distribution (density $f_l(x;\theta^*)$). 
Likewise, $\mu_l(\theta^*)$ \eqref{mapping} can be different from 
$$\theta^*_l=\arg\max_{\theta'_l\in \Theta_l}\int \ln f_l(x;\theta'_l)f_l(x;\theta^*_l)d\lambda(x).$$
In order to see this, note that \eqref{koondumineII2} and \eqref{mapping} in the context of the mixture 
model specialize to
\begin{align}
  \mu_l(\theta,\theta^*)&
   =\arg\max_{\theta'_l\in \Theta_l}\int_{S_l(\theta)}\ln f_l(x;\theta'_l)f(x;\theta^*)d\lambda(x)
\label{eq:mulmix}\\
 \mu_l(\theta) 
&=\arg\max_{\theta'_l\in \Theta_l}\int_{S_l(\theta)}\ln f_l(x;\theta'_l)
                                          \bigl(\sum_i p_if_i(x;\theta_i)\bigr)d\lambda(x),
\label{kass0}
\end{align}
respectively.
We also emphasize that $\Delta$ 
can be significant which justifies the adjustment.
\begin{example}\label{normal} Let
$$f(x;\theta^*)={1\over K}\sum_{l=1}^K\phi(x;\theta^*_l),$$
where $\phi(x;\theta^*_l)$ is the density of the $d$-variate normal distribution with
identity covariance matrix and vector of means $\theta^*_l\in \mathbb{R}^d=\Theta_l$ for 
$l=1,2,\ldots,K$.
In this case, for each $K$-tuple of parameters
$\theta=(\theta_1,\ldots, \theta_K)$, the decision-rule for
the alignment is essentially as follows (disregarding  possible
ties):  $v_{\theta}(x)=i$ if and only if $\|x-\theta_i\|\leq \min_j \|x-\theta_j\|$.
Thus, the decision regions in this case correspond to the Voronoi partition in its original sense,
justifying our generalization of this term.
Now, 
it can be easily seen that for all $m=1,\ldots,d$: 
\begin{equation}\label{koer}
(\mu_l(\theta))_m={\sum_{i=1}^K\int_{S_l(\theta)}x_m\phi(x;\theta_i)dx_1\cdots dx_d \over
\sum_{i=1}^K\int_{S_l(\theta)}\phi(x;\theta_i)dx_1\cdots dx_d}.
\end{equation}
\end{example}
When $d$ and $K$ are large, the exact integration in \eqref{koer} still requires 
intensive computations,  for which reason one may be interested in approximations of \eqref{koer}.
In the context of the above example, one might think of the following
approximations for $\Delta_l(\theta)=\theta_l-\mu_l(\theta)$:
\begin{figure}
\begin{center}
\setlength{\unitlength}{0.00072in}
\begingroup\makeatletter\ifx\SetFigFont\undefined%
\gdef\SetFigFont#1#2#3#4#5{%
  \reset@font\fontsize{#1}{#2pt}%
  \fontfamily{#3}\fontseries{#4}\fontshape{#5}%
  \selectfont}%
\fi\endgroup%
{\renewcommand{\dashlinestretch}{30}
\begin{picture}(7224,6939)(0,-10)
\put(4062,3012){\circle*{150}}
\put(4887,1437){\circle*{150}}
\put(2412,2562){\circle*{150}}
\put(6537,3462){\circle*{150}}
\put(4962,5187){\circle*{150}}
\put(3087,4362){\circle*{150}}
\put(1587,4137){\circle*{150}}
\put(2487,5712){\circle*{150}{150}}
\put(987,3987){\ellipse{150}{150}}
\put(6837,3837){\ellipse{150}{150}}
\put(2937,4512){\ellipse{150}{150}}
\put(4212,3162){\ellipse{150}{150}}
\put(4940,815){\ellipse{150}{150}}
\put(1887,2037){\ellipse{150}{150}}
\put(5037,5637){\ellipse{150}{150}}
\put(2187,6237){\ellipse{150}{150}}
\path(2412,3612)(12,2412)
\path(3012,3312)(3462,1737)
\path(3012,3312)(3462,1737)
\path(3462,1737)(5412,2712)
\path(3462,1737)(5412,2712)
\path(3462,1737)(2712,12)
\path(3462,1737)(2712,12)
\path(5412,2712)(5112,3612)
\path(5412,2712)(5112,3612)
\path(5412,2712)(6912,1737)
\path(5412,2712)(6912,1737)
\drawline(4287,3912)(4287,3912)
\path(4212,3912)(5112,3612)
\path(4212,3912)(5112,3612)
\path(2412,3612)(3012,3312)
\path(2412,3612)(3012,3312)
\path(3012,3312)(4212,3912)
\path(3012,3312)(4212,3912)
\path(5112,3612)(7212,5412)
\path(5112,3612)(7212,5412)
\path(4212,3912)(3687,5262)
\path(4212,3912)(3687,5262)
\path(2412,3612)(2262,4812)
\path(2412,3612)(2262,4812)
\path(2262,4812)(537,5862)
\path(2262,4812)(537,5862)
\path(2262,4812)(3687,5262)
\path(2262,4812)(3687,5262)
\path(3687,5262)(4212,6912)
\path(3687,5262)(4212,6912)
\dashline{60.000}(4962,837)(4287,2187)(4287,2112)
\dashline{60.000}(4737,762)(4437,1362)
\blacken\thicklines
\path(4530.915,1274.793)(4437.000,1362.000)(4450.416,1234.544)(4530.915,1274.793)
\thinlines
\dashline{60.000}(2637,3462)(1887,2037)
\dashline{60.000}(1887,2037)(1362,3087)
\dashline{60.000}(1887,2037)(3312,1362)
\dashline{60.000}(1887,2037)(3237,2412)
\dashline{60.000}(1737,2187)(1512,2637)
\blacken\thicklines
\path(1605.915,2549.793)(1512.000,2637.000)(1525.416,2509.544)(1605.915,2549.793)
\thinlines
\dashline{60.000}(1887,2262)(2037,2562)
\blacken\thicklines
\path(2023.584,2434.544)(2037.000,2562.000)(1943.085,2474.793)(2023.584,2434.544)
\thinlines
\dashline{60.000}(2112,2037)(2637,2187)
\blacken\thicklines
\path(2533.980,2110.765)(2637.000,2187.000)(2509.255,2197.302)(2533.980,2110.765)
\thinlines
\dashline{60.000}(1962,1887)(2187,1737)
\blacken\thicklines
\path(2062.192,1766.122)(2187.000,1737.000)(2112.115,1841.006)(2062.192,1766.122)
\thinlines
\blacken\thicklines
\path(2388.392,5793.863)(2487.000,5712.000)(2466.534,5838.515)(2388.392,5793.863)
\path(2487,5712)(2187,6237)
\blacken\path(2285.608,6155.137)(2187.000,6237.000)(2207.466,6110.485)(2285.608,6155.137)
\thinlines
\blacken\thicklines
\path(1092.503,4059.761)(987.000,3987.000)(1114.331,3972.448)(1092.503,4059.761)
\path(987,3987)(1587,4137)
\blacken\path(1481.497,4064.239)(1587.000,4137.000)(1459.669,4151.552)(1481.497,4064.239)
\put(2637,5637){\makebox(0,0)[lb]{\smash{{{\SetFigFont{12}{14.4}{\rmdefault}{\mddefault}{\updefault}$\theta^*_5$}}}}}
\put(2350,6250){\makebox(0,0)[lb]{\smash{{{\SetFigFont{12}{14.4}{\rmdefault}{\mddefault}{\updefault}$\mu_5$}}}}}
\put(1062,2187){\makebox(0,0)[lb]{\smash{{{\SetFigFont{12}{14.4}{\rmdefault}{\mddefault}{\updefault}$\Delta_1^6$}}}}}
\put(1212,1137){\makebox(0,0)[lb]{\smash{{{\SetFigFont{12}{14.4}{\rmdefault}{\mddefault}{\updefault}$l=1$}}}}}
\put(4287,312){\makebox(0,0)[lb]{\smash{{{\SetFigFont{12}{14.4}{\rmdefault}{\mddefault}{\updefault}$l=2$}}}}}
\put(6612,2862){\makebox(0,0)[lb]{\smash{{{\SetFigFont{12}{14.4}{\rmdefault}{\mddefault}{\updefault}$l=3$}}}}}
\put(5412,4887){\makebox(0,0)[lb]{\smash{{{\SetFigFont{12}{14.4}{\rmdefault}{\mddefault}{\updefault}$l=4$}}}}}
\put(3312,6387){\makebox(0,0)[lb]{\smash{{{\SetFigFont{12}{14.4}{\rmdefault}{\mddefault}{\updefault}$l=5$}}}}}
\put(687,4887){\makebox(0,0)[lb]{\smash{{{\SetFigFont{12}{14.4}{\rmdefault}{\mddefault}{\updefault}$l=6$}}}}}
\put(3312,4662){\makebox(0,0)[lb]{\smash{{{\SetFigFont{12}{14.4}{\rmdefault}{\mddefault}{\updefault}$l=7$}}}}}
\put(4062,3387){\makebox(0,0)[lb]{\smash{{{\SetFigFont{12}{14.4}{\rmdefault}{\mddefault}{\updefault}$l=8$}}}}}
\put(3237,4062){\makebox(0,0)[lb]{\smash{{{\SetFigFont{12}{14.4}{\rmdefault}{\mddefault}{\updefault}$\Delta_7\approx 0$}}}}}
\put(4287,2787){\makebox(0,0)[lb]{\smash{{{\SetFigFont{12}{14.4}{\rmdefault}{\mddefault}{\updefault}$\Delta_8\approx 0$}}}}}
\put(1062,4212){\makebox(0,0)[lb]{\smash{{{\SetFigFont{12}{14.4}{\rmdefault}{\mddefault}{\updefault}$\Delta_6$}}}}}
\put(2050,5800){\makebox(0,0)[lb]{\smash{{{\SetFigFont{12}{14.4}{\rmdefault}{\mddefault}{\updefault}$\Delta_5$}}}}}
\end{picture}
}
\end{center}
\caption{An example of the Voronoi partition for multiple components. True parameters $\theta^*$ and corresponding
$\mu(\theta^*)$ are marked with solid and transparent dots, respectively. For $l=1$ component, 
$\Delta_l^j(\theta^*)$, the correction components, are indicated. For $l=2$ component, the main direction of the 
correction is indicated. 
It appears a reasonable approximation to neglect the corrections for the estimators corresponding to the bounded 
Voronoi regions.}
\label{fig:voronoj}
\end{figure}
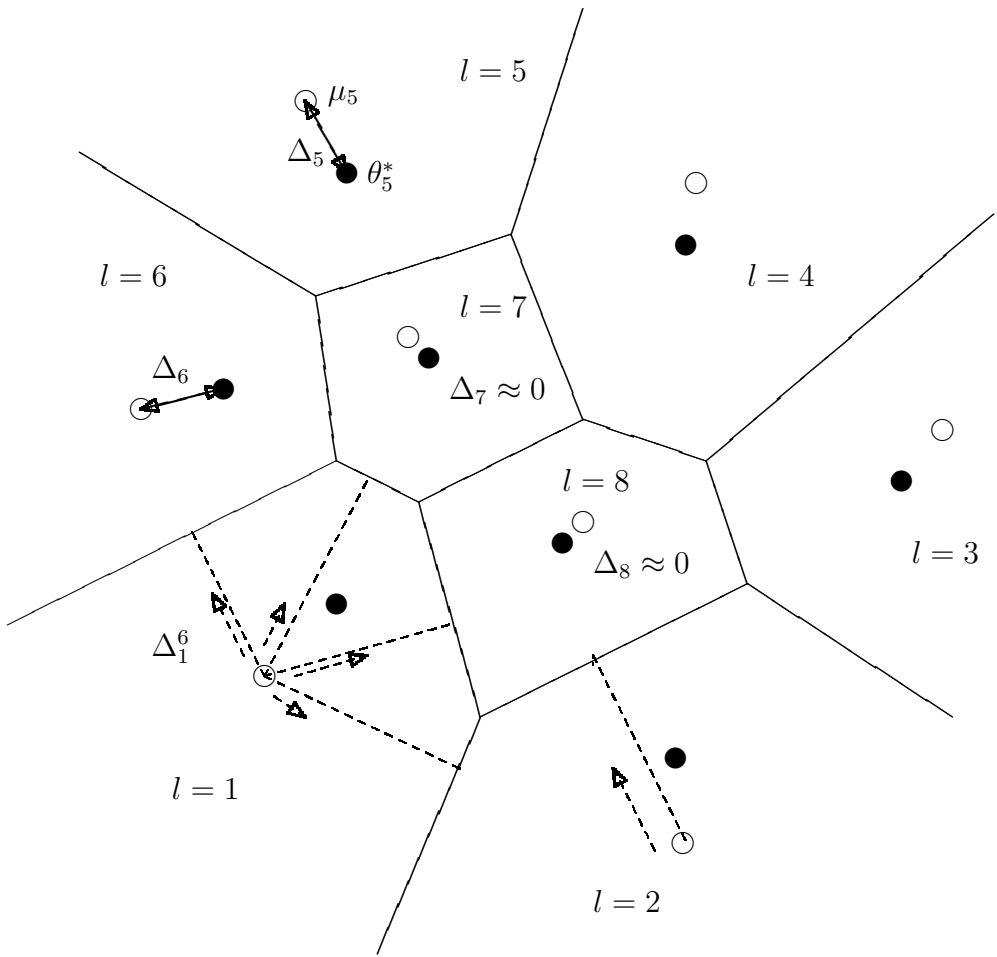
\begin{enumerate}[1.)]
\item Approximate $\Bigr(\sum_l f_i(x;\theta_i)\Bigl)I_{S_l(\theta)}$ in \eqref{koer} by
$f_l(\theta_l,x)I_{S_l(\theta)}$, so
\begin{eqnarray}\label{eq:approxdelta}
(\mu_l)_m&\approx&{\int_{S_l(\theta)}x_m\phi(x;\theta_l)dx_1\cdots dx_d \over
\int_{S_l(\theta)}\phi(x;\theta_l)dx_1\cdots dx_d }. 
\end{eqnarray}
This approximation is based on the limiting case when the components are ``infinitely'' far 
from each other.
   \item If $K>d$, then some components are fully surrounded by others, namely, the partition cells 
corresponding to such ``internal'' components are  bounded (Figure \ref{fig:voronoj}).  It is then
conceivable that  $\Delta_l$'s that  correspond to the bounded cells are less significant than the
others, in which case one might only correct the estimators of the internal components. 
   \item Every Voronoi cell is determined by several hyperplanes and 
every such hyperplane $H^j$ gives rise to $\Delta^j_l$, a component of $\Delta_l$ in the direction 
perpendicular to $H^j$ and corresponding to the $l^{th}$ term in the sum in \eqref{koer}. 
Thus, $\Delta_l=\sum_j\Delta_l^j$ (see, for example, $l=1$ cell in Figure \ref{fig:voronoj}). 
It may be reasonable to find only the "main direction" of
correction, i.e. the largest $\Delta_l^j$ for each $l$ (see, for example, the $l=2$ cell in 
Figure \ref{fig:voronoj}).
\end{enumerate}
\begin{remark}
In Example \ref{normal}, the  decision
regions correspond to the Voronoi partition in its original sense.
Moreover, it is easy to see that in this particular case, the Viterbi
training is none other than the  well-known (generalized) Lloyd algorithm
designed for finding vector quantizers, which in this case are
also called {\em $K$-means} (see, e.g. \citep{loid}). In this
case, the estimators obtained by the Viterbi training are 
empirical $K$-means. The latter estimators enjoy certain desirable properties, and 
in particular they are  consistent
with respect to the population $K$-means \citep{pollard}.  However, they need not 
be consistent  with respect to $\theta^*$, our parameters of interest.  In the
mixture case, the Viterbi training  can always be 
considered as the (generalized) Lloyd algorithm, and the estimators obtained by 
Viterbi training can be regarded as  (generalized) empirical
$K$-means \citep{gray}. This observation links the study of  Viterbi  Training 
and related algorithms to the theory of vector quantization. 
\end{remark}
\section{VA2 -- Adjustment of the second type}
\label{sec:second}
The adjusted Viterbi training is designed  to  asymptotically fix the true parameter $\theta^*$,  
returning approximately the correct solution given this solution as the initial guess  and given 
an infinitely large data sample: VA1$(\theta^*)\approx \theta^*$.  VA2 goes further and attempts to maximally
expand $\{\theta: \text{VA1}(\theta)\approx \theta^*\}$, the set of parameter values that are asymptotically 
mapped to the true ones, to $\{\theta: \text{VA2}(\theta)\approx \theta^*\}$.  Specifically, if the algorithm 
ever  arrives at $\mathcal{S}(\theta^*)$, the Voronoi partition corresponding to the true 
parameters $\theta^*$,  then we would like to coerce the adjusted estimates to return $\theta^*$.   
Let us explain these ideas in more detail.  
\\ \indent
Let ${\mathcal S}^*$ stand for ${\mathcal S}(\theta^*)$, the true Voronoi partition (that also coincides with the 
Bayes decision boundary). The mapping $\theta\mapsto {\mathcal S}(\theta)$  is generally  many-to-one, hence the set 
$\Theta({\mathcal S}^*)=\{\theta: {\mathcal S}(\theta)={\mathcal S}^*\}$ generally contains more than one element.  
(This also means that guessing ${\mathcal S}^*$, i.e. guessing any element from $\Theta({\mathcal S}^*)$, is generally 
easier than guessing $\theta^*$.)  We now introduce VA2:
\\ \indent
Note first that $\mu_l(\theta,\theta^*)$  in \eqref{eq:mulmix}, as well as the estimate $\hat \mu_l(\theta)$ 
in \eqref{hyljes}, depends on $\theta$ through ${\mathcal S}(\theta)$ only.
However, the correction $\Delta_l(\theta)=\theta_l-\mu_l(\theta,\theta)$ does depend on $\theta$ fully and hence
would not generally  work (in the sense of \eqref{eq:approx}) for an arbitrary 
$\theta^j\in \Theta(\mathcal S^*)$ unless $\theta^j=\theta^*$. 
We now attempt to improve the first type of adjustment that is based on adding $\Delta(\theta^j)$ to 
$\hat\mu(\theta^j)$.
Namely, we propose the following iterative update for $l=1,\ldots,K$:
Next, define $\mu_{l,\Theta({\mathcal S}(\theta^0))}(\theta)$ (as function of $\theta$ only) to be the restriction of 
$\mu_l(\theta^0,\theta)$ to $\Theta({\mathcal S}(\theta^0))$, and   write 
$\mu_{l,\theta^0}(\theta)$ in place of the more cumbersome $\mu_{l,\Theta({\mathcal S}(\theta^0))}(\theta)$.  
Let
\begin{eqnarray}\label{sipelgas}
  \theta_l^{j+1}=
\begin{cases}\mu_{l,\theta^j}^{-1}(\hat\mu_l(\theta^j)),
\quad\text{if a unique }\mu_{l,\theta^j}^{-1}(\hat\mu_l(\theta^j))
                                                  \quad\text{exists}\\
              \hat\mu_l(\theta^j)+\Delta_l(\theta^j)\quad\text{otherwise.} \end{cases}
\end{eqnarray}
For any $\theta^j$ and $\theta^*$, the event that  $\hat\mu(\theta^j,\theta^*)$ belongs to 
the range of $\mu(\theta^j,\theta)$ as a function of $\theta\in \Theta({\mathcal S}(\theta^j))$ 
is of zero probability, as Example \ref{normalII} illustrates.  Hence, the introduction of the individual inverses
$\mu_{l,\theta^j}^{-1}$ $l=1,\ldots,K$ is essential, although still not always effective:
In some mixture models (a mixture of normal distributions with unequal weights is 
one such example), 
for a fixed $l$, the event that  $\hat\mu_l(\theta^j)$ belongs in the range of 
$\mu_l(\theta^j,\theta)$ (as function of $\theta\in \Theta({\mathcal S}(\theta^j))$) need not occur with 
probability one for all $\theta^j$ and $\theta^*$.  This, and also the fact that the 
inverses  in general need not have a closed form, or
may require intensive computations, may reduce the attractiveness of the suggested method.  Further
discussion of the computational issues related to this method is outside the scope of this paper, 
except for mentioning the possibility of various, e.g. linear or quadratic, approximations of the 
above functions $\mu_{l,\theta}^{-1}$.   
\\ \indent
In order to better understand the meaning of the new adjustment, imagine 
that $\theta^j\in\Theta({\mathcal S}^*)$.  We would then expect for $l=1,\ldots,K$:
$$\theta^{j+1}_l=\mu_{l,\theta^j}^{-1}(\hat\mu_l(\theta^j))
              =\mu_{l,\theta^*}^{-1}(\hat\mu_l(\theta^*))
              \approx \mu_{l,\theta^*}^{-1}(\mu_l(\theta^*,\theta^*))
              = \theta^*_l.$$
The above argument, of course, also depends on the regularity of the above inverses at 
$\mu_l(\theta^*,\theta^*)$ $l=1,\ldots,K$, and in this regard our experiments 
in \S\ref{sec:simu} provide encouraging results for an important model similar to the model in 
the following example:
\begin{example}\label{normalII}  Let $f(x;\theta^*)={1\over
2}\phi(x-\theta_1^*)+{1\over 2}\phi(x-\theta_2^*)$, where $\phi$ is the
density of the standard normal distribution. In this
case any Voronoi partition is specified by a single parameter $t=0.5(\theta_1+\theta_2)$
solving $\phi(t-\theta_1)=\phi(t-\theta_2)$
(ties are evidently inessential in this context). The true Voronoi partition corresponds to
$t^*=0.5(\theta^*_1+\theta^*_2)$. Given a Voronoi partition ${\mathcal S}(t(\theta))$,   
$\Theta(t)=\{(t-a),(t+a):a\in \mathbb{R}^+\}$.
Hence, restricted to $\Theta(t)$, the function 
$\mu_{{\mathcal S}(t)}(\theta)=(\mu_{1,{\mathcal S}(t)}(\theta),\mu_{2,{\mathcal S}(t)}(\theta))$ depends 
on one parameter only: 
Let $a$ be this parameter and define $\mu_{{\mathcal S}(t)}(\theta(a))=(\mu_1(a),\mu_2(a))$ as follows:
$\mu_1(a)=-a(1-2\Phi(-a))-2\phi(-a)+t$, $\mu_2(a)=2t-\mu_1(a)$,
where $\Phi$ is the distribution function of the standard normal distribution.
After calculating $\hat{\mu}_1<\hat{\mu}_2$ from the data, the inversion equations of \eqref{sipelgas} become
\begin{equation}\label{sipsik}
t-[a(1-2\Phi(-a))+2\phi(a)]=\hat{\mu}_1, \quad  t+[a(1-2\Phi(-a))+2\phi(a)]={\hat\mu}_2.
\end{equation}
Obviously \eqref{sipsik} has a (unique) solution if and only if $\hat{\mu}_1$, ${\hat\mu}_2$ are 
symmetric with respect to $t$ and the probability of this latter event is clearly zero under the model.
Thus, as suggested in \eqref{sipelgas}, we consider the equations separately:
 \begin{align}\label{sipsa}
 a(1-2\Phi(-a))+2\phi(a)&=t-{\hat\mu}_1\\
 \label{sipsb}
 a(1-2\Phi(-a))+2\phi(a)&={\hat\mu}_2-t.
 \end{align}
It can be shown that \eqref{sipsa} and \eqref{sipsb} have  unique solutions, let us denote the latter by 
$a_1$ and $a_2$, respectively.  The points $t-a_1$ and $t+a_2$ will be now taken as the estimators of 
$\theta^*_1$ and $\theta^*_2$ for the next step of iterations.
\end{example}
\begin{center}\underline{VA2}
\end{center}
\begin{enumerate}[1.)]
\item Choose $\theta^0=(\theta^0_1,\ldots,\theta^0_K)$.
\item Given $\theta^j$, find ${\mathcal S}(\theta^j)$
and define empirical  measures $\hat P_l^n(\theta^j)$ as in \eqref{empmix}.
\item For every $\hat P_l^n$, find $\hat \mu_l(\theta^j,\theta^*)$ as in  \eqref{hyljes}.
\item Update $\theta^{j+1}$ in accordance with \eqref{sipelgas}.
\end{enumerate}
\subsection{Unknown weights}
\label{sec:kaalud}
We consider the case when the mixture weights
$p_l$ are unknown, which  corresponds to the case of
the unknown transition parameters $(P,\pi)$ in the general HMM context.
\\ \indent
The Voronoi partition depends on the weight-vector
$p=(p_1,\ldots,p_K)$ as well as on  $\theta$. Hence, ${\mathcal
S}(\theta,p)$ and the vector $p$ should  be reestimated at each
step along with $\theta$. Given a Voronoi partition ${\mathcal
S}=\{S_1,\ldots, S_K\}$, the simplest way to estimate the weights
$p_l$ is to take $p_l=\hat P_n(S_l)$, the empirical measure of $S_l$. Hence all
the algorithms considered so far can be modified accordingly to include the weight estimation 
as in \eqref{phinnang}.
\begin{equation}\label{phinnang}
p^{j+1}_l=\hat P_n(S_l(\theta^j,p^j)),\quad l=1,\ldots, K.
\end{equation}
Taking into account the asymptotics, it is easy to correct the
estimators $p^{j+1}$ as well. Indeed, suppose
$\theta^j=\theta^*$ and $p^j=p$, i.e. ${\mathcal
S}(\theta^j,p^j)={\mathcal S}(\theta^*,p)={\mathcal S}^*$. If $n\to
\infty$, then
\begin{equation}\label{lovi}
\hat P_n\Bigl(S_l(\theta^*,p)\Bigr)\stackrel{\rm a.s.}
{\to}P\Bigr(S_l(\theta^*,p)\Bigl)=\int_{S_l(\theta^*,p)}f(x;\theta^*)d\lambda=
\sum_ip_i\int_{S_l(\theta^*,p)}f_i(x;\theta^*_i)d\lambda.
\end{equation}
In general the latter differs from $p_l$.
The difference is
$p_l-P\bigl(S_l(\theta^*,p)\bigr)$.
Hence, by analogy with \eqref{mapII}, we can define the weight correction 
$D(\theta,p)=(D_1(\theta,p),\ldots,D_K(\theta,p))$ as follows:
\begin{equation}\label{mager}
D_l(\theta,p)=p_l-\sum_ip_i\int_{S_l(\theta,p)}f_i(x;\theta_i)d\lambda,
\end{equation}
which is also data independent.
We now summarize the above by giving a formal definition of the adjusted
 Viterbi training with the weight correction. The
  Viterbi training and the second type of adjusted Viterbi training with $p$  unknown can be defined similarly.
\begin{center}\underline{VA1 with the weight correction}\end{center}
\begin{enumerate}[1.)]
\item Choose $\theta^0=(\theta^0_1,\ldots,\theta^0_K)$ and  $p^0=(p^0_1,\ldots p^0_K)$
\item Given  $\theta^j=(\theta^j_1,\ldots,\theta^j_K)$ and  $p^j=(p^j_1,\ldots p^j_K)$
define  the Voronoi partition ${\mathcal S}(\theta^j,p^j)=\{S_1,\ldots, S_K\}$ as in \eqref{eq:S1} and \eqref{eq:S},
and the empirical measures $\hat P_l^n(\theta^j,p^j)$ as in \eqref{empmix}.
\item Put $\theta^{j+1}=\hat \mu^j(\theta^j)+\Delta(\theta^j),$
where $\hat \mu^j$ is defined in \eqref{kass0}.
\item Put $p^{j+1}=\hat P_n(S_l(\theta^j,p^j))+D(\theta^j,p^j)$.
\end{enumerate}
\section{Simulation studies}\label{sec:simu}
\begin{figure}[htbp]
 \leavevmode 
\begin{center}
\epsfig{file=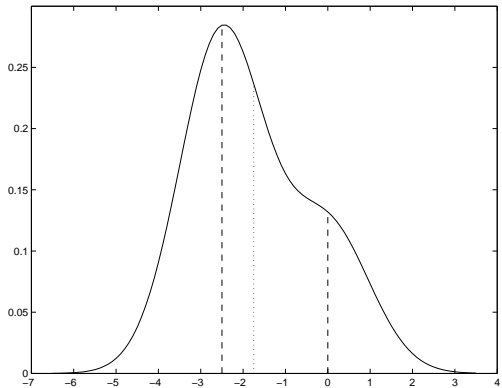, width=2.6in}
\end{center}
\caption{$\frac{1}{\sqrt{2\pi}}(0.7e^{-\frac{(x+2.5)^2}{2}}+0.3e^{-\frac{x^2}{2}})$. 
The dashed vertical lines indicate the means of the individual components, 
the dotted line marks the mean of the mixture.}
\label{fig:mixture}
\end{figure}

In order to support our theory of adjusted Viterbi Training we simulate 
1000 i.i.d. random  samples of size 1000 according to the following
mixture: 
$$\frac{1}{\sqrt{2\pi}}(pe^{-\frac{(x-\theta_1)^2}{2}}+(1-p)e^{-\frac{(x-\theta_2)^2}{2}}).$$
The true parameters in our experiments are $\theta^*=(-2.5,0)$ and  $(p,1-p)=(0.7,0.3)$.
The corresponding density is plotted in Figure \ref{fig:mixture}.
Note that for all such mixtures with $p>0.5$ and $\theta_1<\theta_2$, 
$\theta_2-\theta_1<\sqrt{2 p/(1-p)}$ ($=2.1602$ in our case) implies that the both means 
fall on one side of the decision boundary, which makes detection of the second component particularly
difficult as is already becoming the case in our setting with $\theta^*_2-\theta^*_1$ $=2.5$.
\\ \indent
Our main goal is  to compare the
performances of VT, VA1, and EM in terms of the  accuracy, convergence, 
amount of computations per iteration, and the total amount of computations.
We implement these algorithms in Matlab \citep{matlab},  providing
a fair comparison of their computational intensities based on their  execution times.
Our code is available for the reader's perusal \citep{VA_matlab} and is fully-optimized 
for speed in the case of VT and EM.  Consequently, our simulations possibly only overestimate 
the execution times for VA1.  
\\ \indent
Additionally, we compare VA2 with the above algorithms by the  accuracy and
convergence.  We use a numerical solver to compute the adjustment function of VA2 and 
presently make no effort to replace this by a computationally efficient approximation. 
Hence, we  do not discuss the computational intensity of VA2 in this work.  
\\ \indent
In our experiments, the algorithms are instructed to terminate as soon as the $L_2$ 
distance between consecutive $\theta$ updates falls below $0.001$.  
We also provide a high precision MLE computed with a built-in matlab optimization function.
The cases of known and unknown weights (\S\ref{sec:kaalud}) are considered in 
\S\ref{sec:known_weights} and \S\ref{sec:unknown_weights}, respectively.   
We report the following statistics for each of the algorithms in the form: 
average$\pm$one standard deviation.
\begin{itemize}
\item $\theta=(\theta_1,\theta_2)$ - the estimates of the means;
\item $p$ - the estimate of the weight of the first component;
\item $||\theta-\theta^*||_{1,2}$ -  $L_1$- and $L_2$-normed distances between $\theta$ 
and the true parameters;
\item $n$ - number of steps used by the algorithm;
\item $T$ - total time in milliseconds to execute the entire algorithm;
\item $t$ - time in  milliseconds to execute one iteration of the algorithm;
\end{itemize}
\subsection{Known weights}\label{sec:known_weights}
It is often the case in speech recognition models that the weights  are 
assumed known. 
\\ \indent
First, consider $(-1,2)$ as an ``arbitrary'' initial guess for $\theta$. 
Table~\ref{tab:arbitrary} presents the performance statistics based on the 1000 samples.
The baseline Viterbi method terminates
quickly (on average in $9.04$ steps), outperformed only by VA2, but is the least accurate 
among the considered methods.  
As expected, VT also requires fewest computations: 0.2 ms per iteration and 1.85 ms total.
Ranked from low to high, accuracies of VA1, VA2, and EM appear similar, and are about three times 
superior to that of VT. In units of the VT execution time, EM compares to VA1 as 16.85:6.7
per iteration, and as 20.43:7.59 by the total execution times.
\begin{table}[h]\caption{``Arbitrary'' initial guess.}\label{tab:arbitrary} 
{\footnotesize 
\begin{center}\leavevmode
\begin{tabular}{|c|c|c|c|c|c|}\hline
                           & VT          & VA1         & VA2          & EM           & MEL           \\ \hline 
$\theta_1$  &-2.4869$\pm$0.0497&-2.4952$\pm$0.0500&-2.4959$\pm$0.0498&-2.4970$\pm$0.0456&-2.4973$\pm$0.0456 \\ \hline 
$\theta_2$  &0.2880$\pm$0.0732&0.0099$\pm$0.0917&0.0082$\pm$0.0916&0.0030$\pm$0.0757&0.0024$\pm$0.0757\\ \hline 
$||\theta-\theta^*||_1$&0.3291$\pm$0.0844 &   0.1138$\pm$0.0681    &  0.1133$\pm$0.0678  & 0.0958$\pm$0.0562 & 0.0958$\pm$0.0562 \\ \hline 
$||\theta-\theta^*||_2$&0.2927$\pm$0.0727 &   0.0902$\pm$0.0537    &  0.0899$\pm$0.0536  & 0.0761$\pm$0.0451 & 0.0761$\pm$0.0451 \\ \hline 
$n$  & 9.04$\pm$1.55  &    10.49$\pm$1.61   &   7.84$\pm$1.59  & 11.20$\pm$0.42    & N/A         \\ \hline 
$t$                    & 0.20$\pm$0.05        &     1.34$\pm$0.19    &  39.24$\pm$1.56       &  3.37$\pm$0.07        & N/A         \\ \hline 
$T$  & 1.85$\pm$0.55         &    14.04$\pm$2.95     &  308.57$\pm$68.63        & 37.79$\pm$1.56         & N/A         \\ \hline 
\end{tabular}
\end{center}}
%
\caption{Correct initial guess.} \label{tab:true}
{\footnotesize
\begin{center}\leavevmode
\begin{tabular}{|c|c|c|c|c|c|}\hline
                            & VT          & VA1         & VA2          & EM           & MLE           \\ \hline 
$\theta_1$  &-2.4904$\pm$0.0495&-2.4973$\pm$0.0488&-2.4973$\pm$0.0490&-2.4973$\pm$0.0455&-2.4973$\pm$0.0456 \\ \hline 
$\theta_2$  &0.2820$\pm$0.0729&0.0051$\pm$0.0880&0.0052$\pm$0.0892&0.0024$\pm$0.0753&0.0024$\pm$0.0757\\ \hline 
$||\theta-\theta^*||_1$&0.3223$\pm$0.0829 &   0.1087$\pm$0.0661    &  0.1102$\pm$0.0664  & 0.0953$\pm$0.0561 & 0.0958$\pm$0.0562 \\ \hline 
$||\theta-\theta^*||_2$&0.2867$\pm$0.0721 &   0.0861$\pm$0.0523    &  0.0874$\pm$0.0525  & 0.0756$\pm$0.0450 & 0.0761$\pm$0.0451 \\ \hline 
$n$  & 5.56$\pm$1.72  &     5.06$\pm$1.57   &   4.73$\pm$1.55  &  5.69$\pm$1.29    & N/A         \\ \hline 
$t$                    & 0.22$\pm$0.02        &     1.37$\pm$0.05    &  42.23$\pm$0.94       &  3.42$\pm$0.08        & N/A         \\ \hline 
$T$  & 1.21$\pm$0.31         &     6.91$\pm$2.05     &  199.52$\pm$65.67        & 19.39$\pm$4.28         & N/A         \\ \hline 
\end{tabular}
\end{center}}
\end{table}
In order to illustrate the asymptotic fixed point property, we initialize the algorithms to
$(-2.5,0)$, the true value of the parameters, see Table \ref{tab:true}.
In this case, as expected, the both types of adjustments exit noticeably faster than VT and EM,
are comparable in accuracy to  EM, and are about three times more accurate than VT.  
Unlike VA1, VA2 or EM, the baseline algorithm, as predicted, disturbs the correct initial guess, 
resulting in an appreciable bias.  The times per iteration of VA1 and EM are similar to as
before, and their total times are (in units of the VT time): 5.71 and 16.03, respectively.
\\ \indent
In order to illustrate the idea of the second type of adjustment, we now initialize the
algorithms to $(-3.1229,0.8771)$, which produces the same decision boundary 
$t=-0.9111$ as  $\theta^*=(-2.5,0)$, the true values. Table \ref{tab:true_bound}
collects these results. Note that since VT and VA2 depend on the initial guess only via 
the decision boundary, 
they produce in this case exactly the same results (disregarding a small rounding error) 
as in the case of the correct initial guess (Table \ref{tab:true}).  
As expected, VA2 now terminates significantly faster than its competitors,
and accuracy-wise is  only slightly superior to VA1 and slightly inferior to EM. 
The times per iteration of VA1 and EM  are similar to as before, and their total times 
are 7.84 and 20.83, respectively. 
\subsection{Unknown weights}\label{sec:unknown_weights}
Assume now that the weights are unknown (\S\ref{sec:kaalud}) and hence need to be estimated along 
with the means.  
We use the same data and the same three types of conditions as in the case of known weights:
Arbitrary initialization to $(-1,2)$ (Table \ref{tab:arbitrary_uw}), initialization to the correct
values $(-2.5,0)$ (Table \ref{tab:correct_uw}), and initialization to $(-3.1229,0.8771)$, 
an arbitrary point giving rise the correct intercomponent boundary (Table \ref{tab:correct_bound_uw}).   
VT and the adjusted algorithms VA1 and VA2 in this case are implemented with the asymptotic correction 
\eqref{mager}.
\begin{table}[h]\caption{Correct decision boundary.} \label{tab:true_bound}
{\footnotesize
\begin{center}\leavevmode
\begin{tabular}{|c|c|c|c|c|c|}\hline
                            & VT          & VA1         & VA2          & EM           & MLE           \\ \hline 
$\theta_1$  &-2.4904$\pm$0.0495&-2.4954$\pm$0.0497&-2.4973$\pm$0.0490&-2.4971$\pm$0.0456&-2.4973$\pm$0.0456 \\ \hline 
$\theta_2$  &0.2820$\pm$0.0729&0.0094$\pm$0.0909&0.0052$\pm$0.0892&0.0030$\pm$0.0757&0.0024$\pm$0.0757\\ \hline 
$||\theta-\theta^*||_1$&0.3223$\pm$0.0829 &   0.1131$\pm$0.0668    &  0.1102$\pm$0.0664  & 0.0958$\pm$0.0562 & 0.0958$\pm$0.0562 \\ \hline 
$||\theta-\theta^*||_2$&0.2867$\pm$0.0721 &   0.0897$\pm$0.0528    &  0.0874$\pm$0.0525  & 0.0761$\pm$0.0450 & 0.0761$\pm$0.0451 \\ \hline 
$n$  & 5.56$\pm$1.72  &     7.09$\pm$1.38   &   4.72$\pm$1.56  &  7.44$\pm$0.94    & N/A         \\ \hline 
$t$                    & 0.22$\pm$0.03        &     1.35$\pm$0.05    &  42.37$\pm$1.12       &  3.42$\pm$0.08        & N/A         \\ \hline 
$T$  & 1.22$\pm$0.31         &     9.56$\pm$1.81     &  200.24$\pm$66.30        & 25.41$\pm$3.19         & N/A         \\ \hline 
\end{tabular}
\end{center}}
\caption{Unknown weights.  ``Arbitrary'' guess.} \label{tab:arbitrary_uw}
{\footnotesize
\begin{center}\leavevmode
\begin{tabular}{|c|c|c|c|c|c|}\hline
                            & VT          & VA1         &  VA2          & EM           & MLE           \\ \hline 
$p$  &0.747 $\pm$0.031&0.703 $\pm$0.028&0.702 $\pm$0.028&0.700 $\pm$0.024&0.699 $\pm$0.024 \\ \hline 
$\theta_1$  &-2.4299 $\pm$0.0753&-2.4919 $\pm$0.0596&-2.4930 $\pm$0.0594&-2.4976 $\pm$0.0531&-2.4992 $\pm$0.0532 \\ \hline 
$\theta_2$  &0.3944 $\pm$0.1178&0.0194 $\pm$0.1099&0.0173 $\pm$0.1094&0.0070 $\pm$0.0944&0.0039 $\pm$0.0947\\ \hline 
$||\theta-\theta^*||_1$&0.4775 $\pm$0.1653 &   0.1382 $\pm$0.0851    &  0.1372 $\pm$0.0846  & 0.1179 $\pm$0.0708 & 0.1179 $\pm$0.0710 \\ \hline 
$||\theta-\theta^*||_2$&0.4058 $\pm$0.1237 &   0.1084 $\pm$0.0657    &  0.1076 $\pm$0.0653  & 0.0931 $\pm$0.0558 & 0.0931 $\pm$0.0560 \\ \hline 
$n$  &14.16 $\pm$3.60  &    13.85 $\pm$3.25   &  12.23 $\pm$2.86  & 24.90 $\pm$2.60    & N/A         \\ \hline 
$t$                    & 0.72 $\pm$0.15        &     1.32 $\pm$0.09    &  39.01 $\pm$1.26       &  3.52 $\pm$0.35        & N/A         \\ \hline 
$T$  &10.13 $\pm$3.01         &    18.25 $\pm$4.27     &  478.19 $\pm$117.67        & 87.67 $\pm$12.98         & N/A         \\ \hline  
\end{tabular}
\end{center}}
\end{table}
(The maximization in the high precision MLE is now performed in  the three variables.)
\\ \indent
The adjusted algorithms now converge 1.7 (VA1) and 2 (VA2) times faster than EM, 
and, what is more remarkable, VA1 and VA2 converge even faster than VT.  The per iteration
times of VA1 and EM compare as about 1.8:4.8 for all the initializations, 
and the total times -- as 1.8:8.7 (arbitrary guess),
1.3:6.67 (true values), and 1.73:7.64 (true boundary), all in units of the VT time. 
VA1 and VA2 are again at least three times more accurate than VT in $\theta$ estimation 
and about one standard deviation more accurate than VT in the weight estimation.  
They are also comparable in accuracy to EM.
\begin{table}[h]
\caption{Unknown weights.   Correct guess.} \label{tab:correct_uw}
{\footnotesize \begin{center}\leavevmode
\begin{tabular}{|c|c|c|c|c|c|}\hline
                            & VT          & VA1         &  VA2          & EM           & MLE           \\ \hline 
$p$  &0.737 $\pm$0.030&0.699 $\pm$0.026&0.699 $\pm$0.026&0.699 $\pm$0.023&0.699 $\pm$0.024 \\ \hline 
$\theta_1$  &-2.4526 $\pm$0.0700&-2.4987 $\pm$0.0555&-2.4987 $\pm$0.0557&-2.4991 $\pm$0.0522&-2.4992 $\pm$0.0532 \\ \hline 
$\theta_2$  &0.3537 $\pm$0.1114&0.0058 $\pm$0.1007&0.0060 $\pm$0.1021&0.0038 $\pm$0.0925&0.0039 $\pm$0.0947\\ \hline 
$||\theta-\theta^*||_1$&0.4212 $\pm$0.1467 &   0.1244 $\pm$0.0782    &  0.1263 $\pm$0.0782  & 0.1149 $\pm$0.0701 & 0.1179 $\pm$0.0710 \\ \hline 
$||\theta-\theta^*||_2$&0.3626 $\pm$0.1146 &   0.0978 $\pm$0.0607    &  0.0994 $\pm$0.0607  & 0.0907 $\pm$0.0553 & 0.0931 $\pm$0.0560 \\ \hline 
$n$  & 8.53 $\pm$3.47  &     6.01 $\pm$2.44   &   6.27 $\pm$2.40  & 11.89 $\pm$4.24    & N/A         \\ \hline 
$t$                    & 0.74 $\pm$0.04        &     1.36 $\pm$0.05    &  41.01 $\pm$1.24       &  3.53 $\pm$0.06        & N/A         \\ \hline 
$T$  & 6.27 $\pm$2.42         &     8.13 $\pm$3.19     &  257.35 $\pm$99.70        & 41.84 $\pm$14.81         & N/A         \\ \hline 
\end{tabular}
\end{center}}
\caption{Unknown weights.  Correct boundary.} \label{tab:correct_bound_uw}
{\footnotesize \begin{center}\leavevmode
\begin{tabular}{|c|c|c|c|c|c|}\hline
                            & VT          & VA1         &  VA2          & EM           & MLE           \\ \hline 
$p$  &0.737 $\pm$0.029&0.702 $\pm$0.026&0.700 $\pm$0.026&0.700 $\pm$0.023&0.699 $\pm$0.024 \\ \hline 
$\theta_1$  &-2.4517 $\pm$0.0689&-2.4941 $\pm$0.0573&-2.4972 $\pm$0.0556&-2.4981 $\pm$0.0526&-2.4992 $\pm$0.0532 \\ \hline 
$\theta_2$  &0.3549 $\pm$0.1096&0.0148 $\pm$0.1050&0.0087 $\pm$0.1024&0.0059 $\pm$0.0930&0.0039 $\pm$0.0947\\ \hline 
$||\theta-\theta^*||_1$&0.4218 $\pm$0.1459 &   0.1327 $\pm$0.0779    &  0.1271 $\pm$0.0780  & 0.1164 $\pm$0.0694 & 0.1179 $\pm$0.0710 \\ \hline 
$||\theta-\theta^*||_2$&0.3637 $\pm$0.1132 &   0.1043 $\pm$0.0606    &  0.0999 $\pm$0.0606  & 0.0919 $\pm$0.0548 & 0.0931 $\pm$0.0560 \\ \hline 
$n$  & 8.16 $\pm$3.40  &     7.74 $\pm$2.58   &   6.54 $\pm$2.22  & 12.98 $\pm$4.22    & N/A         \\ \hline 
$t$                    & 0.74 $\pm$0.04        &     1.34 $\pm$0.04    &  41.12 $\pm$1.74       &  3.52 $\pm$0.05        & N/A         \\ \hline 
$T$  & 5.97 $\pm$2.36         &    10.32 $\pm$3.35     &  268.96 $\pm$91.44        & 45.59 $\pm$14.69         & N/A         \\ \hline 
\end{tabular}
\end{center}}
\end{table} 
\subsection{Summary of the results}
\label{sec:results}
VA1 is consistently close in 
accuracy to EM which is always superior to Viterbi Training: 
Specifically, in estimating the means, the gain in accuracy is about three-fold as measured by
$L_1$- and $L_2$-distances, and in estimating the weights, it is about one standard deviation. 
\\ \indent
VA1 always converges almost as fast as VT and noticeably (by 30\% in the case of unknown weights)
faster than EM. 
\\ \indent
When the weights are known, an iteration of VA1 is about six times longer than that
of VT and is more than twice as fast as that of EM.  By  total 
execution,   VA1 is at most eight times slower than  VT and is more than two and a half
times faster than  EM.
\\ \indent
When the weights are unknown, VA1 is at most twice slower than VT and more than two and
a half times faster than EM, per iteration. It is also about 50\% slower than than VT 
and more than four times faster than EM in total times. 
\\ \indent
Accuracy of VA2 is consistently between those of VA1 and EM, and VA2 additionally converges
faster than VA1.
\section{Conclusion}
\label{sec:conclusion}
We have considered the problem of parameter estimation of the emission distribution in Hidden Markov Models in 
connection with the two most common estimation algorithms: the EM algorithm and the Viterbi Training
algorithm.   
We have identified the sources of bias, or lack of consistency, in VT estimation in comparison
with EM (MLE) estimation. In the case of HMM, EM computes 
the MLE, which is often consistent. Trading the EM's accuracy for the VT's ease of computations, one
loses, among other things, the asymptotic fixed point property:  VT no longer holds the true parameter
values fixed,  even asymptotically.  In this work, we have
restored this property, and consequently recovered a certain amount of the EM's accuracy, without a significant 
increase  in computations relative to Viterbi Training.  Specifically, we have derived two types of adjustments
to the baseline Viterbi Training algorithm.  Our first  algorithm, VA1, that we also present as the 
central contribution of this work, is a modification of VT that restores the asymptotic  fixed point 
property. We also present evidence that, at least in the case of mixture models (a special and important 
case of HMM), the price in extra computations for this increase in accuracy can be made reasonable.
The second algorithm, VA2,  in addition to restoring the fixed point property, 
also ensures that the true parameters are asymptotically found  as soon as the true alignment 
(i.e. Voronoi partition) has been found.  
This latter feature may require intensive computations, undermining 
Viterbi Training as a computationally feasible alternative to EM.  We intend to investigate feasible 
approximations to the correction function of VA2 in future work. 
\\ \indent
Certainly, the final decision as to which algorithm to 
use is application dependent, and this work presents valuable information to facilitate 
such selection, especially in the context of the Gaussian mixture models.  For this special case, 
we have provided simulation studies based on 1000 large random samples which illustrate  
the key features of the adjusted algorithms in contrast with EM and baseline Viterbi Training. 
In our simulations, VA1 demonstrates a significant increase of accuracy (three-fold and one standard 
deviation in in estimating the mixture means and weights, respectively) relative to VT.  In fact,
the accuracy of VA1 is already comparable to that of EM.  Computation-wise, VA1 in our studies 
is still several factors faster than EM. We therefore suggest replacing VT by VA1 in applications 
that can afford computing (to a variable precision)  the correction function in appreciation of 
the increased accuracy. 
\bibliographystyle{kluwer}

\end{document}